\documentclass[12pt]{amsart}
\usepackage{amssymb,amsmath,txfonts,amsrefs, amscd}
\usepackage{array,multirow,bigstrut}
\usepackage{appendix,mathtools}
\usepackage[all]{xy}
\usepackage{xargs}  
\usepackage[pdftex,dvipsnames]{xcolor}  
\usepackage[colorinlistoftodos,prependcaption,textsize=tiny]{todonotes}
\usepackage{hyperref} 

\newcommandx{\unsure}[2][1=]{\todo[linecolor=red,backgroundcolor=red!25,bordercolor=red,#1]{#2}}

\newcommandx{\change}[2][1=]{\todo[linecolor=blue,backgroundcolor=blue!25,bordercolor=blue,#1]{#2}}

\newcommandx{\info}[2][1=]{\todo[linecolor=OliveGreen,backgroundcolor=OliveGreen!25,bordercolor=OliveGreen,#1]{#2}}

\newcommandx{\improvement}[2][1=]{\todo[linecolor=Plum,backgroundcolor=Plum!25,bordercolor=Plum,#1]{#2}}
\newcommandx{\thiswillnotshow}[2][1=]{\todo[disable,#1]{#2}}

\newtheorem{thm}{Theorem}[section]
\newtheorem{prop}[thm]{Proposition}
\newtheorem{lemma}[thm]{Lemma}
\newtheorem{remark}[thm]{Remark}

\newtheorem{defi}[thm]{Definition}
\newtheorem{cor}[thm]{Corollary}

\numberwithin{equation}{section}

\begin{document}

\title[]{Rigidity of complete manifolds with weighted Poincar\'e inequality}

\author{Lihan Wang}
\email{lihan.wang@csulb.edu}
\address{Department of Mathematics and Statistics, CSULB}

\date{}


\begin{abstract}
We consider complete Riemannian manifolds which satisfy a weighted Poincar\`e inequality and have the Ricci curvature bounded below in terms of the weight function. When the weight function has a non-zero limit at infinity, the structure of this class of manifolds at infinity are studied and certain splitting result is obtained. Our result can be viewed as an improvement of Li-Wang's result in \cite{LW3}.
\end{abstract}
\smallskip
\maketitle


\section{Introduction}

We are interested the rigidity of complete Riemannian manifolds which satisfy a weighted Poincar\`e inequality and have the Ricci curvature bounded below in terms of the weight function. P. Li and J. Wang initiated the rigidity study of this class of manifolds in \cite{LW3} and their results recover the rigidity results for manifolds with positive spectrum in \cite{LW1} and \cite{LW2}. Let us recall the following definitions.

\begin{defi}
Let $M^n$ be an $n-$dimensional complete Riemannian manifold. We say that $M^n$ satisfies a weighted Poincar\`e inequality with a nonnegative weight function $\rho(x)$, if the inequality 
\[
\int_{M} \rho(x) \phi^2(x) dV \leq \int_{M} |\nabla \phi|^2 dV
\] is valid for all compactly supported smooth function $\phi$. 

Moreover, if the $\rho-$metric, defined as 
\[
ds^2_{\rho}=\rho ds^2_M
\]is complete, we say that $M^n$ has the property $(\mathcal{P}_{\rho})$.
\end{defi}

In \cite{LW3}, P. Li and J. Wang proved manifolds satisfying the property $(\mathcal{P}_{\rho})$ do not split when the weight function $\rho$ goes to zero at infinity :
\begin{thm}[P. Li \& J. Wang]\label{thm7.2}
Let $M^n$ be a complete manifold of dimension $n\geq 4$ with the property $\mathcal{P}_{\rho}$. Suppose the Ricci curvature of $M^n$ satisfies the lower bound
\[
\rm{Ric}_{M}\geq -\frac{4}{n-1}\rho.
\] If the weight function $\rho$ satisfies that 
\[
\underset{x\rightarrow \infty}{\lim}\rho(x)=0,
\] then $M$ has only one end.
\end{thm}
 A natural question to ask is how many ends the manifold will have when $\underset{x\rightarrow \infty}{\lim}\rho(x) \neq 0$. The main purpose of this paper is to answer this question. At first, we have the following observation  about warped product metrics.  
\begin{prop}\label{prop}
Let $M^n=R\times N^{n-1}$ be a complete Riemannian manifold with the warped product metric 
\[
ds_M^2=dt^2+\eta ^2\left( t\right) ds_N^2.
\] Here $N^{n-1}$ is compact and $\eta > 0$. Assume that $M^n$ satisfies the property $(P_\rho)$ with $\rho=\frac{\left( n-1\right) ^2}4\eta ^{-1}\eta ^{\prime \prime }
$, and has the Ricci curvature bounded below:
\[
Ric_M \geq -\frac 4{n-1}\rho .
\]

Let $E$ be the end with infinity volume and $F$ be the end with finite volume. Then
\[
\limsup_{ F \cap \infty} \rho \left( x\right) \geq
\liminf_{ E\cap \infty} \rho \left( x\right) >0.
\]In particular, if there exist limits as
\[
\lim_{ E\cap \infty} \rho \left( x\right) =A,\;\;\lim_{ 
F\cap \infty} \rho \left( x\right) =B,
\]then $A\leq B$.
\end{prop} 

In view of this proposition, it becomes of special interest to prove a
splitting theorem when $\rho $ has a nonzero limit at infinity on a Riemannian manifold $M$. The main result we get is as follows:

\begin{thm}[Main theorem]\label{mthm}
Let $M^n$ be a complete Riemannian manifold of dimension $n\geq 4$ with the property $(P_\rho)$. Assume the weight function $\rho $ has non-zero limit at infinity and 
\[
Ric_M \geq -\frac 4{n-1}\rho.
\]
Then $M^n$ has
\begin{enumerate}
\item either only one end, in particular which has infinite volume;
\item or one end with infinite volume and one end with finite volume. In this case, $M^n$ is given by $M^n=\mathbb{R}\times N^{n-1}$ with the warped product metric 
 \[
ds_M^2=dt^2+\eta ^2\left( t\right) ds_N^2,\] for some positive function $\eta(t)$, and some compact manifold $N^{n-1}$. Moreover, $\rho$ is a function of $t$ alone satisfying 
\[
\frac{\left( n-1\right) ^2}4\eta ^{-1}\eta ^{\prime \prime }=\rho .
\]
\end{enumerate}
\end{thm}
The case of dimension $3$ has been studied by P. Li and J. Wang in \cite{LW3}. In fact, they proved that $M^3$ either has only one end, or splits with two ends with infinite volume, or splits with one end with finite volume and one end with infinite volume. See Theorem B in \cite{LW3} for details.

We follow the spirit in \cite{LW3} to prove our main theorem. The strategy is to use harmonic functions constructed in \cite{LT} for manifolds with more than one end. The splitting is from reading the equality related to harmonic functions. When $\rho$ goes to zero at infinity, boundary terms will converge to zero on each end in \cite{LW3} which imply the desired equality. However this does not happen when $\rho$ has a non-zero limit. We employ different analysis and find out that involved estimates will have different signs on different ends. Clever manipulations with careful choices of cut-off functions finally derive the desired equality.  

We want to point out that parabolicity and nonprabolicity are used to distinguish ends in \cite{LW3}. In fact, under our assumption, one end has finite volume if and only if it is parabolic. One end has infinite volume if and   only if it is nonparabolic. See Section 2 for details.

The paper is organized as follows. We will recall some concepts and results from \cite{LW3} in Section 2. Then we prove the Main Theorem in Section 3 and the Proposition \ref{prop} in Section 4.

\begin{remark}
P. Li and J. Wang also studied Riemannian manifolds with property $(\mathcal{P}_{\rho})$ and $\rm{Ric} \geq -\frac{n-1}{n-2} \rho$ in \cite{LW3}. Their result shows that $M^n$ either has only one non-parabolic end or split as a warped product with two non-parabolic ends when $\rho$ has a non-zero limit at infinity. 
\end{remark}

{\bf Acknowledgements.} This paper is dedicated to Professor Peter Li for the occasion of his seventieth birthday. The author is very grateful for his
valuable guidance, scholarly input and consistent encouragement. The author would like to thank Ovidiu Munteanu for suggesting this question, sharing ideas and helpful discussions.

\section{Preliminary} 
In this section, we recall some concepts and known results for the proof of  Theorem \ref{mthm} and Proposition \ref{prop}. More details can be found in \cite{LW3}. 
\subsection{$\rho-$ distance}
Consider a Riemannian manifold $M^n$ which satisfies the property $\mathcal{P}_{\rho}$. With respect to the $\rho-$ metric 
\[ds_{\rho}^2=\rho ds_M^2,\]  define the $\rho-$ distance function:
\[
r_{\rho}(x,y)=\underset{\gamma}{\inf} \,l_{\rho}(\gamma),
\] which is the infimum of the length of all smooth curves joining $x$ and $y$ with respect to $ds^2_{\rho}$. 
It is easy to see that 
\[|\nabla r_{\rho}|^2(x)=\rho(x).\] 
Throughout this paper, we denote the ball with respect to the $\rho-$ distance by \[
B_{\rho}(p, R)=\{x\in M^n |r_{\rho}(p, x)<R \}
\] comparing to the usual ball
\[B(p, R)=\{x\in M^n |r(p, x)<R\}.\]
When $p\in M^n$ is a fixed point, we will suppress the dependency of $p$ and write $B_{\rho}(R)=B_{\rho}(p, R)$ and $B(R)=B(p, R)$. 
 If $E$ is an end of $M^n$, we denote $E_{\rho}(R)=B_{\rho}(R)\cap E$.
We also define that 
\[
S(R)=\underset{B_{\rho}(R)}{\sup}\,\sqrt{\rho}
\]to be the supremum of $\sqrt{\rho}$ over the set $B_{\rho}(R)$.
\subsection{About ends}
Recall an end is simply an unbounded component of $M$ outside of a compact smooth domain of $M^n$. Say one end $E$ is {\it nonparabolic} if there exists a symmetric, positive, Green's function for the Laplacian acting on $L^2$ functions with Neumann boundary condition on $\partial E$. Otherwise, we say that $E$ is a {\it parabolic} end. In \cite{LW3}, P. Li and J. Wang proved the following geometric condition for parabolicity and nonparabolicity of ends (see Corollary 3.2 in \cite{LW3}).

\begin{prop}\label{cor32}
Let $E$ be an end of a complete Reimannian manifold $M^n$ with Property $(P_{\rho})$. If $E$ is nonparabolic, then it must have at least quadratic volume growth. In particular, if the weight function $\rho$ satisfies $\underset{x\rightarrow\infty }{\lim \inf} \rho(x)>0$ on $E$, then $E$ is nonparabolic if and only if $E$ has infinite volume.
\end{prop}
Therefore, under our assumptions of Property $(P_{\rho})$ and $\underset{x\rightarrow \infty}{\lim}\rho >0$,  one end $E$ has infinite (finite) volume if and only if it is nonparabolic (parabolic). 

According to this proposition, Corollary 2.3 in \cite{LW3}, and results in \cite{LT}, we have following properties about harmonic functions on ends.
\begin{lemma}[Harmonic functions on ends] \label{hf} Let $M^n$ be a complete Riemannian manifold satisfying property $(\mathcal{P}_{\rho})$. And $\rho$ has the nonzero limit at infinity. Assume $M^n$ has one end with infinite volume and one end $F$ with finite volume. Let $E=M^n\setminus F$. Then there exists a positive harmonic function $f$ such that 
\begin{enumerate}
\item $0<f<1$ on $E$ with $\underset{E}{\inf} f=0$. And there exists a constant $C$ such that
\begin{align*}
&\int_{E_{\rho}(R+1)\backslash E_{\rho}(R)}\rho f^2\leq C \exp(-2R)\\
&\int_{E_{\rho}(R+1)\backslash E_{\rho}(R)}|\nabla f|^2 \leq C \exp(-2R).
\end{align*}
\item $f$ is proper on $F$ with $\underset{ F\cap \infty}{\liminf} f = \infty$. 
\end{enumerate}
\end{lemma}

%
\subsection{About harmonic functions}
We summarize some basic properties of harmonic functions on Riemannian manifolds. The first one is about level sets of harmonic functions, Lemma 5.1 in \cite{LW3}.
\begin{lemma}[Level sets of harmonic functions]\label{lev}
Let $M^n$ be a complete manifold. Let $f >0$ be a bounded harmonic function. Let us denote the level set of $f$ at $t$ by 
\begin{align*}
l(t)&=\{x \in M^n | f(x)=t\}
\end{align*} for $\inf f < t< \sup f$ and we denote the set 
\begin{align*}
\mathcal{L}(a, b)&=\{x\in M^n | a< f(x)< b\}
\end{align*} for $\inf f<a <b< \sup f$. Then 
\begin{align*}
&\int_{\mathcal{L}(a, b)}|\nabla f|^2=(b-a)\int_{l(b)}|\nabla f|.
\end{align*}Moreover,
\begin{align*}
&\int_{l(b)}|\nabla f|=\int_{l(t)}|\nabla f|
\end{align*} for all $\inf f <t< \sup f$.
\end{lemma}

The next one is the improved Bochner formula for harmonic functions, Lemma 4.1 in \cite{LW3}. It implies the splitting when the equality is realized.
\begin{lemma}[Improved Bochner formula]\label{lem4.1}
Let $M^n$ be a complete Riemannian manifold of dimension $n\geq 2$. Assume that the Ricci curvature of $M^n$ satisfies the lower bound\[\rm{Ric}(x)\geq -(n-1)\tau(x)\] for all $x\in M^n$. Suppose $f$ is a non-constant harmonic function defined on $M^n$. Then the function $|\nabla f|$ must satisfy the differential inequality
\[
\Delta |\nabla f| \geq -(n-1)\tau |\nabla f|+\frac{|\nabla |\nabla f||^2}{(n-1)|\nabla f|}
\] in the weak sense. Moreover, if equality holds, then $M^n$ is given by $M^n=\mathbb{R}\times N^{n-1}$ with the warped product metric \[ds^2_M=dt^2+\eta^2(t)ds^2_N\] for some positive function $\eta$. In this case, $\tau(t)$ is a function of $t$ alone satisfying \[\eta^{\prime\prime}\eta^{-1}(t)=\tau(t).\]
In addition, $N^{n-1}$ must be compact if $M^n$ has more than one end.
\end{lemma}

Let us recall the local gradient estimate of Cheng-Yau \cite{CY75}(see \cite{LW2}) with respect to our assumption on Ricci curvature. 

\begin{lemma}[Local gradient estimate]\label{gra1}

Let $M^n$ be a complete Riemannian manifold with $\rm{Ric}\geq -\frac{4}{n-1}\rho(x)$. For any positive harmonic function $f(x)$, there is
\[
|\nabla f|(x) \leq \left(2\underset{B(x, r)}{\sup}\sqrt{\rho}+Cr^{-1}\right)f(x)
\] where $C$ is a constant depending only on $n$. 
\end{lemma}
This local gradient estimate can be stated in terms of $S(R)=\underset{B_{\rho}(R)}{\sup}\sqrt{\rho}$.
 \begin{cor}\label{gra2} Let $M^n$ be a complete Riemannian manifold with $\rm{Ric}\geq -\frac{4}{n-1}\rho(x)$. For any positive harmonic function $f$ and $x\in B_{\rho}(R)$, there is 
  \begin{align*}
|\nabla \ln f(x)| &\leq CS(R+1)
\end{align*} where $C$ is a constant depending only on $n$.
 \end{cor}
 \begin{proof}
Let $F(r)=r-\left(\underset{B(x, r)}{\sup} \sqrt{\rho}\right)^{-1}$.  Since $F$ becomes negative as $r \rightarrow 0$ and $F$ goes to infinity as $r \rightarrow \infty$, there there exists $R_0$ such that $F(R_0)=0$, i.e.\[R_0=\left(\underset{B(x, R_0)}{\sup} \sqrt{\rho}\right)^{-1}.\] Moreover, $B(x, R_0)\subset B_{\rho}(x,1)$. Let $x \in B_{\rho}(R)$. It follows that $B(x, R_0) \subset B_{\rho}(x, 1)\subset B_{\rho}(R+1)$. Then choose $r=R_0$ in Lemma \ref{gra1} and we get 
\begin{align*}
|\nabla \ln f(x)|&\leq \left(2 \underset{B_(x,R)}{\sup}\sqrt{\rho}+CR_0^{-1}\right) \leq \left(2 \underset{B_{\rho}(R+1)}{\sup}\sqrt{\rho}+CR_0^{-1}\right) \\
&=(2+ C) \underset{B_{\rho}(R+1)}{\sup}\sqrt{\rho}=CS(R+1)
\end{align*} for some constant $C$ which depends only on $n$.
\end{proof}
\subsection{The rigidity theorem in \cite{LW3}}
In the end, we cite Theorem 5.2 in \cite{LW3} which we will use to exclude the case of two ends with infinite volume , i.e., two nonparabolic ends, in our main theorem. 
\begin{thm}\label{thm5.2}
Let $M^n$ be a complete manifold with dimension $n\geq 3$. Assume that $M$ satisfies property $(\mathcal{P}_{\rho})$ for some nonzero weight function $\rho\geq 0$. Suppose the Ricci curvature of $M$ satisfies the lower bound 
\[
\rm{Ric} (x)\geq -\frac{n-1}{n-2}\rho(x)
\]for all $x \in M$. If $\rho$ satisfies the growth estimate 
\[
\underset{R\rightarrow \infty}{\lim}\inf \frac{S(R)}{F(R)}=0
\]where 
\begin{equation*}
F(R)=\left\{ 
\begin{array}{c}
\exp(\frac{n-3}{n-2}R) \\ 
R 
\end{array}%
\begin{array}{l}
\text{when }n\geq 4 \\ 
\text{when }n=3\
\end{array}%
\right.
\end{equation*}
then either 
\begin{enumerate}
\item $M$ has only one non-parabolic end; or
\item M has two non-parabolic ends and is given by $M=\mathbb{R}\times N$ with the warped product metric 
\[
ds^2_M=dt^2+\eta^2(t)ds^2_N,
\] for some positive function $\eta(t)$, and some compact manifold $N$. Moreover, $\rho(t)$ is a function of $t$ alone satisfying 
\[
(n-2)\eta^{\prime\prime}\eta^{-1}=\rho.
\]
\end{enumerate}

\end{thm}


\section{Proof of Main Theorem}

For convenience, we write 
$\underset{x\rightarrow \infty}{\lim}\rho \left( x\right) =A^2$ for some constant $A>0$ according to the assumption on $\rho$. Then it implies that \[\underset{R\rightarrow \infty}{\lim}S(R)=A.\] Therefore the growth assumption in Theorem \ref{thm5.2} is fulfilled. We also notice that $-\frac 4{n-1}\rho > -\frac{n-1}{n-2}\rho$  for $n\geq 4$. Then Theorem \ref{thm5.2} asserts that $M^n$ must have at least one nonparabolic end, i.e. one end with infinite volume. 
\begin{remark}In fact, according to Corollary 1.4 in \cite{LW3}, the assumption of the weighted Poincar\'e inequality alone is enough to imply that the manifold has at least one nonparabolic end. 
\end{remark}

 At first, let us assume that $M^n$ has two ends with infinite volume, i.e., two nonparabolic ends. Theorem \ref{thm5.2} 
 asserts that $M^n$ must be given the warped product metric \[
ds^2_M=dt^2+\eta^2(t)ds^2_N\] and $
(n-2)\eta^{\prime\prime}\eta^{-1}=\rho$. However this implies that $\rm{Ric}=-\frac{n-1}{n-2}\rho$ along $\partial t$ direction which is impossible for $n\geq 4$. So $M^n$ has only one end of infinite volume when $n\geq 4$.

Hence we may assume that $M^n$ has one end of infinite volume and one end of finite volume. Let $F$ denote the end of finite volume and $E=M^n\setminus F$. Then according to Lemma \ref{hf}, there will exist a positive harmonic function $f$ satisfying: 
 \begin{enumerate}
\item $0<f<1$ on $E$ with $\underset{E}{\inf} f=0$. And there exists a constant $C$ such that
\begin{align}
&\int_{E_{\rho}(R+1)\backslash E_{\rho}(R)}\rho f^2\leq C \exp(-2R)\label{decay1}\\
&\int_{E_{\rho}(R+1)\backslash E_{\rho}(R)}|\nabla f|^2 \leq C \exp(-2R).\label{decay2}
\end{align}
\item $f$ is proper on $F$ with $\underset{ F\cap \infty}{\liminf} f = \infty$. 
\end{enumerate}
 
 Let $g=|\nabla f|^{\frac{1}{2}}$. Lemma \ref{lem4.1} implies that 
\begin{align}\label{ine1}
\Delta g&\geq -\frac{2}{n-1}\rho g-\frac{n-3}{n-1}g^{-1}|\nabla g|^2.
\end{align} We will show that this inequality is an equality. Then Lemma \ref{lem4.1} will assert that $M^n$ splits and is given by $M^n=\mathbb{R}\times N^{n-1}$ with $
ds_M^2=dt^2+\eta^2(t) ds_N^2$, and $-\frac{4}{n-1}\rho=\rm{Ric}_{11}=-(n-1)\eta^{\prime\prime}$ in the $\frac{\partial}{\partial t}$ direction.

Let 
\[h=\Delta g +\frac{2}{n-1}\rho g+\frac{n-3}{n-1}g^{-1}|\nabla g|^2.\] Here $h \geq 0$. We only need to show $h=0$ in order to show the equality in \eqref{ine1}.

Let $\phi$ be some cut-off function chosen later. With the integration by parts, it follows:
\begin{align*}
\int_M\phi^2gh&=\int_M\phi^2 g\left(\Delta g +\frac{2}{n-1}\rho g+\frac{n-3}{n-1}g^{-1}|\nabla g|^2\right)\\
&=\frac{-2}{n-1}\int_M\phi^2|\nabla g|^2+\frac{2}{n-1}\int_M\rho\phi^2g^2-2\int_M\phi g\langle \nabla \phi, \nabla g\rangle.
\end{align*}
 Apply the weighted Poincar\'e inequality to the term $\int_M \rho\phi^2g^2$. Then we get
\begin{align}\label{eq2}
0\leq \int_M\phi^2gh
&\leq \frac{2}{n-1}\int_M|\nabla\phi|^2g^2-\frac{2(n-3)}{n-1}\int_M \phi g\langle \nabla \phi, \nabla g\rangle.
\end{align} We will get $h=0$ if we can prove that the right hand side of \eqref{eq2} is nonpositive as the cut off function $\phi$ approaches to $1$, i.e., 
\begin{align}\label{ine2}
\mathcal{H}=\int_M\left| \nabla \phi \right| ^2g^2-\left( n-3\right) \int_M\phi g\langle\nabla
\phi,\nabla g\rangle\leq 0\;as\;\;\phi \rightarrow 1. 
\end{align}

To prove \eqref{ine2}, we choose the cut-off function $\phi =\chi \psi$ with $\chi =\varphi ^{\frac 32}$ as follows: 
\[
\varphi =\left\{ 
\begin{array}{c}
0\;\;\;\;\;\;\;\;\;\;on\;\;L(0,\delta \varepsilon )\cap E ,\\ 
\frac{\log f-\log \delta \varepsilon }{-\log \delta }\;\;on\;\;L(\delta
\varepsilon ,\varepsilon )\cap E, \\ 
1\;\;\;\;\;\;\;\;\;on\;L(\varepsilon ,\infty )\cap E,
\end{array}
\right. 
\] 

\[
\varphi =\left\{ 
\begin{array}{c}
0\;\;\;\;\;\;\;\;\;\;on\;\;L(\beta T,\infty )\cap F ,\\ 
\frac{\log \beta T-\log f}{\log \beta}\;\;on\;\;L(T,\beta T)\cap F, \\ 
1\;\;\;\;\;\;\;\;\;on\;L(0,T)\cap F,
\end{array}
\right. 
\]
and
\[
\psi =\left\{ 
\begin{array}{c}
1\;\;on\;\;E_\rho (R-1)\cup F, \\ 
R-r_{\rho}\;\;on\;E_\rho (R)\backslash E_\rho (R-1), \\ 
0\;\;on\;E\backslash E_\rho \left( R\right)
\end{array}
\right. 
\] with $0< \delta< \epsilon<1$ and $1<T< \beta$. Here $\psi=1$ on $F$ since the level set of $f$ on F is always compact due to properness of $f$ on $F$. We would like to point it out that $\chi =\varphi ^{\frac 32}$ is essential for our argument when $\rho$ has non-zero limit at infinity.

Define \[\Omega=\{x\in M| \phi(x)\neq 0, \nabla \phi \neq 0\}.\] With above choice of $\phi$, we get:
\begin{align*}
\mathcal{H}
&=\int_{\Omega}\left| \nabla \phi \right| ^2g^2-\frac{n-3}{2} \int _\Omega\phi\chi \langle\nabla \psi,\nabla g^2\rangle+\frac{n-3}{4}\int _\Omega(\psi^2 g^2\Delta \chi^2+g^2\langle\nabla \psi^2, \nabla \chi^2\rangle)\\
&-\frac{n-3}{4}\int_{\partial \Omega}\chi^2_n\psi^2g^2.
\end{align*}
Since $|\nabla \phi|^2=\chi^2|\nabla \psi|^2+\psi^2|\nabla \chi|^2+\frac{1}{2}\langle\nabla \chi^2,\nabla\psi^2\rangle \geq \frac{1}{2}\langle\nabla \chi^2,\nabla\psi^2\rangle$, it follows that
\begin{align}\label{ine3}
\begin{aligned}
\mathcal{H}&\le \frac{n-1}{2}\int_{\Omega}\left| \nabla \phi \right| ^2g^2-\frac{n-3}{2} \int _\Omega\phi\chi \langle\nabla \psi, \nabla g^2\rangle\\
&+\frac{n-3}{4}\int _\Omega \psi^2 g^2\Delta \chi^2-\frac{n-3}{4}\int_{\partial \Omega}\chi^2_n\psi^2g^2.
\end{aligned}
\end{align}
We are going to prove that the first, second and fourth term in the right hand of \eqref{ine3} go to zero and the third term will become nonpositive as the cut off function $\phi$ approaches to $1$. 

Let us list two facts frequently used in following proof at first. According to Corollary \ref{gra2} and the assumption $\underset{x\rightarrow \infty}{\lim} \rho=A^2 >0$, for $x\in \Omega 
$, there is
\begin{align*}
|\nabla \ln f|(x)&\leq CS(R+1)\rightarrow CA\, \, \text{as $R\rightarrow \infty$}.
\end{align*}
Another fact is that $\int_{l(t)}|\nabla f|$ is a constant independent of $t$ by Lemma \ref{lev}.

\begin{lemma} The first term in \eqref{ine3} goes to zero:
\[
\int_{\Omega}\left| \nabla \phi \right| ^2g^2\rightarrow 0,\;\;as\;\;\phi \rightarrow 1. 
\]
\end{lemma}


\begin{proof}

 
 By the choice of $\phi$, it is easy to see that 
\begin{align}\label{ine4}
\int_\Omega \left| \nabla \phi \right| ^2g^2 \leq \frac{9}{2}\left((\ln\delta)^{-2}\int_{E\cap \Omega}+(\ln\beta)^{-2}\int_{F\cap \Omega}\right)|\nabla \ln f|^2 |\nabla f|+2\int_{E\cap \Omega}\rho |\nabla f|.
\end{align}
 Using the co-area formula, we get
\begin{align*}
& \int_{\Omega\cap E}|\nabla \ln f|^2 |\nabla f| \leq CS(R+1)\int_{\Omega\cap E}|\nabla \ln f| |\nabla f|\leq CS(R+1)(-\ln \delta)\int_{l(t)}|\nabla  f|,\\
& \int_{\Omega\cap F}|\nabla \ln f|^2 |\nabla f| \leq CS(R+1)\int_{\Omega\cap F}|\nabla \ln f| |\nabla f|\leq CS(R+1)(\ln \beta)\int_{l(t)}|\nabla  f|.
 \end{align*} 
Then it follows that \begin{align*}
&\frac{9}{2}\left((\ln\delta)^{-2}\int_{E\cap \Omega}+(\ln\beta)^{-2}\int_{F\cap \Omega}\right)|\nabla \ln f|^2 |\nabla f|\\
& \leq \frac{9}{2}CS(R+1)\left((-\ln \delta)^{-1} +(\ln \beta)^{-1} \right)\int_{l(t)}|\nabla  f|\rightarrow 0\, \, \text{as $\phi \rightarrow 1$.}
\end{align*} 
Here we use the fact that $\int_{l(t)}|\nabla  f|$ is a constant  independent of $t$. Thus the first part in the right hand side of \eqref{ine4} goes to zero.

For the second part of the right hand side of \eqref{ine4}, we apply decay properties of $f$ and $\nabla f$ on $E$ by Lemma \ref{hf}. In fact, there is
 \begin{align*}
 \int_{E\cap \Omega} \rho |\nabla f|&\leq \left(\int_{E\cap\Omega} |\nabla f|^2\right)^{1/2}\left(\int_{E\cap\Omega} \rho^2 \right)^{1/2} \\
 &\leq C\exp(-R)S(R)\left(\int_{E\cap\Omega} \rho \right)^{1/2} \; \;\text{by \ref{decay2},}\\
&\leq C\exp(-R) S(R)\left((\delta \varepsilon)^{-2} \int_{\Omega \cap E} \rho f^2\right)^{1/2}\\
& \leq C (\delta\varepsilon)^{-1}S(R)\exp(-2R) \; \;\text{by \ref{decay1}, }\\
&\rightarrow 0 \; \text{as $R \rightarrow \infty$ with $\delta, \varepsilon$ fixed.  }
\end{align*} Thus $\int_{\Omega}|\nabla \phi|^2 g^2 \rightarrow 0$ as $\phi \rightarrow 1$.

\end{proof}

\begin{lemma} The second term in \eqref{ine3} goes to zero:
\[
\int _\Omega\phi\chi \langle\nabla \psi, \nabla g^2 \rangle\rightarrow 0\;\; \text{as $\phi \rightarrow 1$}. 
\]
\end{lemma}
\begin{proof}
By the choice of $\psi, \chi$ and $\phi$, we notice
\begin{align}\label{ine6}
\begin{aligned}
\int _\Omega\phi\chi \langle\nabla \psi, \nabla g^2\rangle&\leq \int _{E\cap \Omega}\sqrt{\rho}|\nabla g^2|\leq \left(\int_{E\cap \Omega} \rho\right)^{1/2}\left(\int_{E\cap \Omega}|\nabla g^2|^2\right)^{1/2}\\
&\leq (\delta \varepsilon)^{-1}\left(\int_{E\cap \Omega} \rho f^2\right)^{1/2}\left(\int_{E\cap \Omega}|\nabla g^2|^2\right)^{1/2}\\
&\leq (\delta \varepsilon)^{-1}\exp(-R)\left(\int_{E\cap \Omega}|\nabla g^2|^2\right)^{1/2} \; \; \, \text{by \ref{decay1} }.
\end{aligned}
\end{align} 
We claim that 
\begin{align}\label{ine7}
\int_{E_{\rho}(R+1)/E_{\rho}(R-2)}|\nabla g^2|^2& \leq  CS^2(R+1)\exp(-2R).
\end{align}
Then \eqref{ine6} and \eqref{ine7} imply that 
\begin{align*}
\int _{\Omega}\phi\chi \langle\nabla \psi,\nabla g^2\rangle&\leq  CS(R+1)\exp(-2R) \left(\delta \varepsilon\right)^{-1} \rightarrow 0\\
& \text{as $ R \rightarrow \infty$ with $\delta, \varepsilon$ fixed.}
\end{align*} Thus, $\int _{\Omega}\phi\chi \langle\nabla \psi,\nabla g^2\rangle \rightarrow 0$ as $\phi \rightarrow 1$.

 To prove \eqref{ine7}, we notice that \eqref{ine1} implies that 
 \begin{align}\label{cla}
 \Delta g^2&\geq -\frac{4}{n-1}\rho g^2.
 \end{align} Choose the cut-off function  $\tau$ as 
\[
\tau=\left\{ 
\begin{array}{c}
r_{\rho}-R+2\;\;\;\;on\; E_{\rho}(R-1)\backslash E_{\rho}(R-2), \\ 
1 \;\;\;\;\;on\;\;E_{\rho}(R)\backslash E_{\rho}(R-1), \\ 
R-r_{\rho}+1\;\;on\;E_{\rho}(R+1)\backslash E_{\rho}(R), \\ 
0\;\;\;\;otherwise.
\end{array}
\right. 
\]  Then multiply \eqref{cla} with $\tau^2g^2$ and run the integration by parts:
\begin{align*}
-\frac{4}{n-1}\int_{M}\tau^2 \rho g^4&\leq \int_{M} \tau^2 g^2 \Delta g^2\\
&=-\int_{M}2g^2 \tau \langle\nabla \tau, \nabla g^2\rangle-\int_{M}\tau^2|\nabla (g^2)|^2\\
&\leq 2 \int_M |\nabla \tau|^2g^4-\frac{1}{2}\int_M \tau^2 |\nabla (g^2)|^2. 
\end{align*} Hence 
\begin{align*} 
\int_M \tau^2 |\nabla g^2|^2& \leq \frac{8}{n-1}\int_M \tau^2 \rho g^4+4\int_M |\nabla \tau|^2g^4.
\end{align*} It follows 
\begin{align*}
 \int_{E_{\rho}(R)\backslash E_{\rho}(R-1)}  |\nabla g^2|^2& \leq \left(\frac{8}{n-1}+4\right)\int_{E_{\rho}(R+1)\backslash E_{\rho}(R-2)}\rho g^4\\
&\leq CS^2(R+1)\int_{E_{\rho}(R+1)\backslash E_{\rho}(R-2)} |\nabla f|^2\\
 &\leq  CS^2(R+1)\exp(-2R), \; \text{by \eqref{decay2}.}
\end{align*}
 
\end{proof}

\begin{lemma} The fourth term in \eqref{ine3} goes to zero:
\[
 \int_{\partial \Omega}\chi^2_n\psi^2g^2\rightarrow 0\;\;as\;\;\phi \rightarrow 1. 
\]
\end{lemma}
\begin{proof}
We notice that $\chi_n|_{l(\varepsilon)}=\frac{3}{2}\varphi\left(\frac{|\nabla \ln f|}{-\ln \delta}\right)$ and $\chi_n|_{l(T)}=-\frac{3}{2}\varphi\left(\frac{|\nabla \ln f|}{\ln \beta}\right)$. In addition, by Corollary \ref{gra2}, $|\nabla \ln f|_{\partial (E_{\rho}(R)\setminus E_{\rho}(R-1) )} \leq  CS(R+1)$. Then it follows that
\begin{align*}
|\int _{\partial \Omega}(\chi^2)_n\psi^2g^2|&\leq 2\int_{\partial \Omega} |\chi_n||\nabla f|\\
&\leq CS(R+1)\left(\int_{l(\epsilon)}\frac{|\nabla f|}{-\ln \delta}+\int_{l(T)}\frac{|\nabla f|}{\ln \beta}\right)\\
&= CS(R+1)\left((\ln \beta)^{-1}+(-\ln \delta)^{-1}\right)\int_{l(t)}|\nabla f| \rightarrow 0, \, \text{as $\phi \rightarrow 1$}
\end{align*} since $\int_{l(\epsilon)}|\nabla f| =\int_{l(T)}|\nabla f| =\int_{l(t)}|\nabla f|$ is a constant independent of $t$.
\end{proof}

\begin{lemma}The third term in \eqref{ine3} becomes nonpositive: 
\[\int_{\Omega}\psi^2g^2\Delta \chi ^2 \leq 0\, \, \text{as $\phi \rightarrow 1$} .\]
\end{lemma} 
\begin{proof}
Different from other terms discussed above, both of $\int_{\Omega\cap E}g^2\Delta \chi ^2$ and $\int_{\Omega\cap F}g^2\Delta \chi ^2$ do not go to zero as $\phi \rightarrow 1$ due to the assumption on $\rho$. We will do different analysis here. In fact it is enough to prove that
\begin{align}\label{ine8}
&\limsup_{\beta ,T\rightarrow \infty }\int_{L\left( T,\beta T\right) \cap
F}g^2\Delta \chi ^2+\limsup_{\varepsilon ,\delta \rightarrow 0}\int_{L\left(
\delta \varepsilon ,\varepsilon \right)\cap{E} }g^2\Delta \chi ^2\leq 0. 
\end{align}
The direct calculation yields
\begin{align*}
\Delta \chi ^2
&=\left\{
\begin{array}{l}
3 |\nabla \ln f|^2\left(-\varphi^2(-\ln \delta)^{-1}+2(-\ln \delta)^{-2}\varphi\right)<0,\; \text{on $E\cap \Omega$}\\
3 |\nabla \ln f|^2\left(\varphi^2(\ln \beta)^{-1}+2(\ln \beta)^{-2}\varphi\right)>0,\;\text{on $F\cap \Omega$}.
\end{array}
\right.
\end{align*} Thus \eqref{ine8} is equivalent to
\begin{equation}\label{ine9}
\begin{aligned}
&\limsup_{\beta ,T\rightarrow \infty }\int_{L(T, \beta T)\cap F}|\nabla f|^3  f^{-2}\left(\varphi^2(\ln \beta)^{-1}+2(\ln \beta)^{-2}\varphi\right)\\
&+\limsup_{\varepsilon ,\delta \rightarrow 0}\int_{L\left(
\delta \varepsilon ,\varepsilon \right)\cap{E} }|\nabla f|^3  f^{-2}\left(-\varphi^2(-\ln \delta)^{-1}+2(-\ln \delta)^{-2}\varphi\right)\leq 0.
\end{aligned}
\end{equation}
With the co-area formula and Corollary \ref{gra2}, it follows that 
\begin{align*}
\int_{L(T, \beta T)\cap F} (\ln \beta)^{-2}\varphi|\nabla f|^3f^{-2}
&\leq  (\ln \beta)^{-2}CS(R+1)\int_{L(T, \beta T)} |\nabla f|^2f^{-1}\\
&\leq CS(R+1)(\ln \beta)^{-1} \int_{l(t)}|\nabla f|\rightarrow 0\, \,\text{as $\beta \rightarrow \infty$},\\
\int_{L(\delta \epsilon, \epsilon)\cap E}(-\ln \delta)^{-2} \varphi|\nabla f|^3f^{-2}&\leq (-\ln \delta)^{-2}CS(R+1)\int_{L(\delta \epsilon, \epsilon)\cap E} \varphi|\nabla f|^2f^{-1}\\
&\leq CS(R+1)(-\ln \delta)^{-1} \int_{l(t)}|\nabla f|\rightarrow 0\,\, \text{as $\delta \rightarrow 0$}.
\end{align*}Here we use the fact that $\int_{l(t)}|\nabla f|$ is a constant independent of $t$. 
Therefore we can deduce the inequality \eqref{ine9} if we can prove
that 
\begin{align}\label{ine10}
\limsup_{\beta ,T\rightarrow \infty }(\ln \beta)^{-1}\int_{L\left(
T,\beta T\right) \cap F}\varphi ^2f^{-2}\left| \nabla f\right| ^3&\leq
\liminf_{\varepsilon ,\delta \rightarrow 0}(-\ln \delta)^{-1}
\int_{L\left( \delta \varepsilon ,\varepsilon \right) }\varphi
^2f^{-2}\left| \nabla f\right| ^3.  
\end{align}

In order to prove \eqref{ine10}, we will derive the upper bound of its left hand side and the lower bound of its right hand side. Recall the gradient estimate in Lemma \ref{gra1} :
\[\left| \nabla \ln f\right| \left( x\right) \leq  2\sup_{B(x,R)}\sqrt{\rho }+CR^{-1}\; \text{for any $R>0$.}\]  Since $\underset{x\rightarrow\infty}{\lim} \sqrt{\rho}=A$, this gradient estimates implies that $\forall \theta >0$,$%
\,\exists \,T_\theta >0$ such that for any  $T>T_\theta $, there is 
\begin{align}\label{grad4}
&\left| \nabla \ln f\right| \left( x\right) \leq 2A+\theta ,\;\;\;x\in L\left( T,\infty \right) \cap F. 
\end{align}
Apply \eqref{grad4} and the co-area formula to the left hand side of (\ref{ine10}). We get the following estimate: 
\begin{align}\label{ine11}
(\ln\beta)^{-1}\int_{L\left( T,\beta T\right) \cap F}\varphi
^2f^{-2}\left| \nabla f\right| ^3 &\leq \frac{2A+\theta }{\ln \beta }%
\int_{L\left( T,\beta T\right) }\varphi^2 f^{-1}\left| \nabla f\right|^2 \\
&\leq\frac{2A+\theta }{3}\int_{l(t)}|\nabla f|.
\end{align}

Next we will estimate the right hand side of (\ref{ine10}) from below. Consider another cut-off function
\[
\tilde{\varphi}=\left\{ 
\begin{array}{c}
0\;\;\;\;on\;L\left( 0,\delta \varepsilon \right) \cap E, \\ 
\varphi \;\;\;\;\;on\;\;L(\delta \varepsilon ,\varepsilon )\cap E, \\ 
\frac{\ln 2\varepsilon -\ln f}{\ln 2}\;\;on\;L(\varepsilon ,2\varepsilon
)\cap E, \\ 
0\;\;\;\;otherwise.
\end{array}
\right. 
\]
Since $\tilde{\varphi}\psi $ has the compact support in $M$, we can apply the weighted Poincare inequality to get the following: 
\begin{align}\label{ine12}
\begin{aligned}
\int_{L\left( \delta \varepsilon ,\varepsilon \right) }\rho \varphi ^2\psi^{2}f
&\leq \int_M\rho \left( \tilde{\varphi}\psi f^{\frac 12}\right) ^2\leq
\int_M\left| \nabla (\tilde{\varphi}\psi f^{\frac 12})\right| ^2\\
&=\frac14\int_{L\left( \delta \varepsilon ,2\varepsilon \right) }\tilde{\varphi}%
^2\psi ^2\frac{\left| \nabla f\right| ^2}f+ 
\int_{L\left( \delta \varepsilon ,2\varepsilon \right) }\tilde{\varphi}%
\psi \langle\nabla f,\nabla (\tilde{\varphi}\psi )\rangle+\int_{L\left( \delta
\varepsilon ,2\varepsilon \right) }\left| \nabla (\tilde{\varphi}\psi
)\right| ^2f \\
&\leq \frac 14\int_{L\left( \delta \varepsilon ,\varepsilon \right)
}\varphi ^2\frac{\left| \nabla f\right| ^2}f+M\left( \delta ,\varepsilon
,R\right)
\end{aligned}
\end{align}where 
\begin{align*}
M\left( \delta ,\varepsilon ,R\right) &=\frac 14\int_{L\left( \varepsilon
,2\varepsilon \right) }\tilde{\varphi}^2\psi ^2\frac{\left| \nabla f\right|
^2}f+\int_{L\left( \delta \varepsilon ,2\varepsilon \right) }\tilde{\varphi}%
^2\psi \langle \nabla f, \nabla \psi \rangle +\int_{L\left( \delta \varepsilon
,2\varepsilon \right) }\tilde{\varphi}\psi ^2 \langle\nabla f,\nabla \tilde{%
\varphi}\rangle \\
&+2\int_{L\left( \delta \varepsilon ,2\varepsilon \right) } \psi^2\left| \nabla 
\tilde{\varphi}\right| ^2f+2\int_{L\left( \delta \varepsilon ,2\varepsilon
\right) }\tilde{\varphi}^2\left| \nabla \psi \right| ^2f.
\end{align*}

{\it Claim:} When $\delta$ and $\varepsilon$ are fixed,
\begin{align}\label{claim}
&\lim_{R\rightarrow \infty }M\left( \delta ,\varepsilon ,R\right) \leq C, 
\end{align}
where $C$ is a constant not depending on $\varepsilon$ and $\delta$. 

Indeed, we have the following estimates on $M\left( \delta ,\varepsilon ,R\right)$ using the co-area formula and decay properties \eqref{decay1} and \eqref{decay2}:
\begin{align*}
\int_{L\left( \varepsilon ,2\varepsilon \right) }\tilde{\varphi}^2\psi ^2%
\frac{\left| \nabla f\right| ^2}f& \leq \int_{L\left( \varepsilon
,2\varepsilon \right) }\frac{\left| \nabla f\right| ^2}f=(\ln 2)\int_{l(t)}|\nabla f|, \\
\int_{L\left( \delta \varepsilon ,2\varepsilon \right) }\tilde{\varphi}\psi
^2\langle\nabla f, \nabla \tilde{\varphi}\rangle& \leq \int_{L\left( \delta
\varepsilon ,2\varepsilon \right) }\left| \nabla f\right| \cdot \left|
\nabla \tilde{\varphi}\right| \\
&=\frac 1{-\ln \delta }\int_{L\left( \delta \varepsilon ,\varepsilon
\right) }\frac{\left| \nabla f\right| ^2}f+\frac 1{\ln 2}\int_{L\left(
\varepsilon ,2\varepsilon \right) }\frac{\left| \nabla f\right| ^2}f=2\int_{l(t)}|\nabla f|,
\\
\int_{L\left( \delta \varepsilon ,2\varepsilon \right) }\tilde{\varphi}%
^2\psi \langle\nabla f,\nabla \psi \rangle
&\leq \int_{L\left( \delta \varepsilon
,2\varepsilon \right) \cap (E_\rho (R)\backslash E_\rho (R-1))}\sqrt{\rho}\left| \nabla
f\right| \\
&\leq \left(\int_{L\left( \delta \varepsilon
,2\varepsilon \right) \cap (E_\rho (R)\backslash E_\rho (R-1))}\rho\right)^{1/2}\left( \int_{L\left( \delta \varepsilon
,2\varepsilon \right) \cap (E_\rho (R)\backslash E_\rho (R-1))}\left| \nabla
f\right|^2\right)^{1/2}\\
&\leq C(\delta \varepsilon)^{-1}\left(\int_{L\left( \delta \varepsilon
,2\varepsilon \right) \cap (E_\rho (R)\backslash E_\rho (R-1))}\rho f^{2}\right)^{1/2}\exp(-R)\\
&\leq C(\delta \varepsilon)^{-1}\exp \left( -2R\right) \rightarrow 0,\;\text{as $R\rightarrow \infty$}
\end{align*} and 
\begin{align*}
\int_{L\left( \delta \varepsilon ,2\varepsilon \right) } \psi^2\left| \nabla 
\tilde{\varphi}\right| ^2f &\leq (-\ln \delta)^{-2}\int_{L\left( \delta \varepsilon ,\varepsilon \right)\cap E(R)\setminus E(R-1) } \frac{\left| \nabla  f\right|^2}{f}+ (\ln2)^{-2} \int_{L\left(  \varepsilon ,2\varepsilon \right) } \frac{\left| \nabla  f\right|^2}{f}\\
&\leq (-\ln \delta)^{-2}\left(\int_{L\left( \delta \varepsilon ,\varepsilon \right)\cap E(R)\setminus E(R-1) }|\nabla f|^2 \right)^{1/2}\left(\int_{L\left( \delta \varepsilon ,\varepsilon \right)\cap E(R)\setminus E(R-1)  }\frac{|\nabla f|^2}{f^2} \right)^{1/2}\\
&+(\ln 2)^{-1}\int_{l(t)}|\nabla f|,\\
&\leq (-\ln \delta)^{-2}C\exp(-R)(\delta \varepsilon)^{-1/2}\left(\int_{L\left( \delta \varepsilon ,\varepsilon \right) }\frac{|\nabla f|^2}{f} \right)^{1/2}+(\ln 2)^{-1}\int_{l(t)}|\nabla f|,\\
&\leq (-\ln \delta)^{-3/2}C\exp(-R)(\delta \varepsilon)^{-1/2}(\int_{l(t)}|\nabla f|)^{1/2}+(\ln 2)^{-1}\int_{l(t)}|\nabla f|\\
&\rightarrow (\ln 2)^{-1}\int_{l(t)}|\nabla f|, \; \text{as $R\rightarrow \infty$}\\
\int_{L\left( \delta \varepsilon ,2\varepsilon
\right) }\tilde{\varphi}^2\left| \nabla \psi \right| ^2f &\leq\int_{L\left( \delta \varepsilon ,2\varepsilon
\right)\cap E_{\rho}(R)\setminus E_{\rho}(R-1) }\rho f\\
&\leq (\delta \varepsilon)^{-1}\left(\int_{L\left( \delta \varepsilon ,2\varepsilon
\right)\cap E_{\rho}(R)\setminus E_{\rho}(R-1) }\rho f^2\right)\\
& \leq C(\delta \varepsilon)^{-1}\exp(-2R)\rightarrow 0, \; \text{as $R\rightarrow \infty$}.
\end{align*}In above estimates, we notice that $\int_{l(t)}|\nabla f|$ is a constant independent of $t$. Therefore our claim is true.

Since $\psi\rightarrow1$ as $ R\rightarrow \infty$, \eqref{claim} together with \eqref{ine12} implies that
\begin{align}\label{ineq12}
\begin{aligned}
\frac 1{-\ln \delta }\int_{L\left( \delta \varepsilon ,\varepsilon \right)
}\rho \varphi ^2f &\leq \frac 14\frac 1{(-\ln \delta) }\int_{L\left( \delta
\varepsilon ,\varepsilon \right) }\varphi ^2\frac{\left| \nabla f\right| ^2}%
f+\frac C{-\ln\delta } \\
&= \frac{\int_{l(t)}|\nabla f|}{12}+\frac C{-\ln \delta},\; \; \text{by the co-area formula}.
\end{aligned}
\end{align}
Using the assumption that $\underset{x\rightarrow \infty}{\lim} \rho=A^2$, we can conclude:
$\forall \theta >0,$ $\exists \varepsilon _\theta >0$ such that if 
$\varepsilon <\varepsilon _\theta $, then \eqref{ineq12} becomes
\begin{align}\label{ineq13}
&\frac{1}{-\ln \delta }\int_{L\left( \delta \varepsilon ,\varepsilon \right)
}\varphi ^2f\leq \frac{\int_{l(t)}|\nabla f|}{3\left( 2A-\theta \right) ^2}+\frac{C}{(-\ln \delta)(A-\theta/2)^2}. 
\end{align}

The next step is to use \eqref{ineq13} and the Schwarz inequality to estimate the right hand side
of (\ref{ine10}) from below: for any $\delta \varepsilon <\varepsilon _0<\varepsilon$, it follows that 
\begin{align*}
\frac{C_2}{3}&=\frac 13\int_{l\left( \varepsilon _0\right) }\left| \nabla f\right| =\frac
1{-\ln \delta }\int_{L\left( \delta \varepsilon ,\varepsilon \right)
}\varphi ^2\frac{\left| \nabla f\right| ^2}f \\
&\leq \left( \frac 1{-\ln \delta }\int_{L\left( \delta \varepsilon
,\varepsilon \right) }\varphi ^2\frac{\left| \nabla f\right| ^3}{f^2}\right)
^{\frac 23}\left( \frac 1{-\ln \delta }\int_{L\left( \delta \varepsilon
,\varepsilon )\right) }\varphi ^2f\right) ^{\frac 13} \\
&\leq \left( \frac 1{-\ln \delta }\int_{L\left( \delta \varepsilon
,\varepsilon \right) }\varphi ^2\frac{\left| \nabla f\right| ^3}{f^2}\right)
^{\frac 23}\left( \frac{\int_{l(t)}|\nabla f|}{3\left( 2A-\theta \right) ^2} +\frac{C}{(-\ln \delta)(A-\theta/2)^2}\right) ^{\frac 13},
\end{align*}

which proves that 
\begin{align}\label{ine13}
\underset{\delta \rightarrow 0}{\liminf}(-\ln \delta)^{-1}\int_{L\left( \delta \varepsilon ,\varepsilon \right)
}\varphi ^2\frac{\left| \nabla f\right| ^3}{f^2}\geq \frac{2A-\theta }%
3\int_{l(t)}|\nabla f| . 
\end{align}

As a consequence, since $\theta $ was arbitrary, (\ref{ine13}) and (\ref{ine11}) show that (\ref{ine10})
holds. 
\end{proof}

By Lemma 3.2, 3.3, 3.4 and 3.5, we can conclude that the right hand side of  \eqref{ine3} goes to nonnegative as the cut off function $\phi \rightarrow 1$. This implies the inequality \eqref{ine2}, i.e., $\mathcal{H} \leq 0$ as $\phi \rightarrow 1$. Then it follows that $0\leq \int_{M}\phi^2gh \leq 0$ in \eqref{eq2} as $\phi \rightarrow 1$. Therefore $h=0$, i.e.,\[\Delta g +\frac{2}{n-1}\rho g+\frac{n-3}{n-1}g^{-1}|\nabla g|^2=0.\] So the inequality \eqref{ine1} becomes an equality. Then according to the argument below \eqref{ine1}, we can conclude the splitting result stated in the theorem.

\section{Warped product : proof of Proposition \ref{prop}}
%
%
%
%
%

We need the following lemma to prove Proposition \ref{prop}.
\begin{lemma}\label{lema}
Suppose a $C^1$ function $y=y(t)>0$ on $[t_0, \infty)$ with $t_0>0$ satisfies 
\begin{align}\label{iney}
&y^{\prime}+y^2 \geq a^2\; \text{on $[t_0, \infty)$}
\end{align} for some $a>0$. Then for any $\epsilon >0$, there exists $t_{\epsilon}>0$ such that 
\[
y\geq a-\epsilon \;\text{on $[t_{\epsilon}, \infty)$}
\]
\end{lemma}
\begin{proof}
Suppose that there exists an interval $(\alpha ,\beta )\subset \;[t_0,\infty )$
such that $0<y<a$ on $(\alpha ,\beta)$. Then this implies that $y^{\prime}>0$ on $(\alpha ,\beta)$ according to \eqref{iney}. In addition, integrating \eqref{iney} on $(\alpha ,\beta)$ yields $y(\beta) > y(\alpha)$. Then $y(\alpha) <a$ since $y^{\prime}>0$ and $0<y<a$ on $(\alpha ,\beta)$. Therefore $0<y<a$ on $[\alpha, \beta)$. Continuing the argument it follows that $y<a$ on $[t_0,\beta)$. 

Then we conclude that either there exists $t_1>t_0$ such
that $y\geq a$ on $[t_1,\infty
) $, or such $t_1$ does not exist which means $0<y<a$ on $[t_0,\infty)$.  Assume the second case is true. Then according to \eqref{iney}, for any $t>t_0$, it follows that
\begin{align*}
\int^{t}_{t_0}\frac{y^{\prime}}{a-y}\geq \int^{t}_{t_0}(a+y)> a(t-t_0).
\end{align*} This implies that there
\[
y(t)\geq a +e^{a(t_0-t)}(a-y(t_0))>a 
\] which is a contradiction to the assumption that $0<y<a$ on $[t_0,\infty)$. Therefore there exists $t_1>t_0$ such that $y\geq a$ on $[t_1,\infty)$. 
\end{proof}

With the similar argument, we can get the following lemma.
\begin{lemma}\label{lemb}
Suppose a $C^1$ function $y=y(t)>0$ on $(- \infty, -t_0]$ satisfies 
\begin{align}\label{ineyb}
&y^{\prime}+y^2 \leq b^2\; \text{on $(- \infty, -t_0)$}
\end{align} for some $b>0$. Then for any $\epsilon >0$, there exists $t_{\epsilon}>0$ such that 
\[
y\leq b+\epsilon \;\text{on $(- \infty, -t_{\epsilon})$}. 
\] 
\end{lemma}

 Now we are ready to prove the Proposition \ref{prop}.\\
\textbf{Proof.} Suppose $E=(0,\infty )\times N$ and $F=(-\infty ,0)\times N.$ According to Theorem 6.3 in \cite{LW3}, our assumption here implies that $\underset{E\cap\infty}{\liminf} \rho>0$ with similar argument. 

 Let 
\[
g\left( t\right) =\eta ^{-\frac{n-1}2}\left( t\right) , 
\]
then $g$ satisfies 
\[
\Delta g=-\frac 2{n-1}\rho g-\frac{n-3}{n-1}g^{-1}\left| \nabla g\right| ^2. 
\]
Consider the cut-off $\phi :M\rightarrow R,$%
\[
\phi \left( t\right) =\left\{ 
\begin{array}{c}
0\;\;\;on\;\left( -\infty ,2L\right) , \\ 
\frac{2L-t}L\;\;on\;\;[2L,L], \\ 
1\;\;on\;\;(L,T), \\ 
\frac{2T-t}T\;\;on\;\;[T,2T], \\ 
0\;\;on\;\;(2T,\infty ),
\end{array}
\right. 
\]
where $L<0$ and $T>0$ are fixed.
Then 
\[\int_M\phi^2g(\Delta g+\frac 2{n-1}\rho g+\frac{n-3}{n-1}g^{-1}\left| \nabla g\right| ^2)=0.
\]
Run the integration by parts and apply the weighted Poincare inequality. It follows that 
\begin{align*}
0&=\int_M\left(-\phi^2|\nabla g|^2-2\phi g\langle\nabla\phi, \nabla g\rangle+\frac 2{n-1}\rho \phi^2g^2+\frac{n-3}{n-1}\phi^2\left| \nabla g\right| ^2\right)\\
& \leq \int_M\left(-\phi^2|\nabla g|^2-2\phi g \langle\nabla\phi, \nabla g\rangle+\frac 2{n-1}|\nabla (\phi g)|^2+\frac{n-3}{n-1}\phi^2\left| \nabla g\right| ^2\right)\\
& = \int_M\left(\frac{2}{n-1}|\nabla \phi|^2g^2-\frac{2(n-3)}{n-1}\phi g\langle \nabla \phi, \nabla g\rangle \right).
\end{align*}That is
\begin{align}\label{ine14}
&\int_M\left| \nabla \phi \right| ^2g^2-\left( n-3\right) \int_M\phi g\langle\nabla
\phi , \nabla g\rangle\geq 0. 
\end{align}

First we claim that
\begin{align}\label{cla1}
&\int_M\left| \nabla \phi \right| ^2g^2\rightarrow 0\;\;\;as\;\;T\rightarrow
\infty ,\;\;L\rightarrow -\infty.  
\end{align}
According to Propositioin \ref{cor32}, the end E is nonparabolic since $\underset{E\cap\infty}{\liminf} \rho>0$. Then a theorem of Varopoulos \cite{V} (see \cite{LW3}) asserts that 
\[
\int^{\infty}_1 A^{-1}(r) dr < \infty
\] where $A(r)$ is the area of the $\partial B(r)$. This implies that 
\begin{align}\label{ine42}
&\int^{\infty}_{-\infty}\eta^{-(n-1)} dt <\infty.
\end{align} Then \eqref{cla1} is true by direct calculations.

Therefore from \eqref{ine14}, it follows
\begin{equation}\label{ine15}
\underset{ T\rightarrow \infty\,  
L\rightarrow -\infty} {\limsup} \int_M\phi g \langle\nabla \phi, \nabla g\rangle\leq 0.
\end{equation}
And the cross term can be written as
\begin{align}\label{ine46}
\begin{aligned}
\int_M\phi g\langle\nabla \phi, \nabla g\rangle &=\frac 12\int_M\phi\langle\nabla \phi, \nabla g^2\rangle \\
&=\frac{n-1}{2T}\int_T^{2T}\frac{2T-t}T\frac{\eta ^{\prime }}\eta dt-\frac{%
n-1}{-2L}\int_{2L}^L\frac{2L-t}L\frac{\eta ^{\prime }}\eta dt.
\end{aligned}
\end{align}

Let us investigate $\frac{\eta^{\prime}}{\eta}$. We notice that $\eta >0$, and $\eta ^{\prime \prime }>0$ implied by $\rho >0$. Then $\eta ^{\prime}$ is increasing. In addition, there is $\int^0_{-\infty}\eta^{-(n-1)} dt<\infty$ by \eqref{ine42}. It follows that $\eta^{\prime}>0$ for $t>0$.

Let us denote $\underset{E\cap\infty }{\liminf}\, \rho=\mathcal{A}>0$ and $\underset{ F \cap\infty }{\limsup}\,\rho=\mathcal{B}$. If $\mathcal{B}=\infty$, then the Proposition follows. Let us assume both of $\mathcal{A}$ and $\mathcal{B}$ are finite. Then for any $\delta>0$, there exists $t_{\delta}$ such that $\eta ^{-1}\eta ^{\prime \prime }\geq \frac{4}{(n-1)^2}\mathcal{A}-\delta>0$ on $%
[t_{\delta},\infty )$, and $\eta ^{-1}\eta ^{\prime \prime }\leq \frac{4}{(n-1)^2}\mathcal{B}+\delta$ on $%
(\infty, -t_{\delta}]$. 

We notice that $y=\frac{\eta^{\prime}}{\eta}$ satisfies that $y^{\prime}+y^2=\frac{\eta^{\prime\prime}}{\eta}=\frac{4}{(n-1)^2}\rho$. Then according to Lemma \ref{lema} and Lemma \ref{lemb}, for any $\varepsilon>0$, there exists $t_{\varepsilon}>t_{\delta}$ such that 
\begin{align}\label{ine48}
\begin{aligned}
\frac{\eta^{\prime}}{\eta}&\geq \sqrt{\frac{4}{(n-1)^2}\mathcal{A}-\delta}-\varepsilon \; \text{on $[t_{\varepsilon}, \infty)$},\\
\frac{\eta^{\prime}}{\eta}&\leq \sqrt{\frac{4}{(n-1)^2}\mathcal{B}+\delta}+\varepsilon \; \text{on $(-\infty, -t_{\varepsilon}]$}.
\end{aligned}
\end{align} 

Now, apply \eqref{ine46}, \eqref{ine48} to \eqref{ine15}. It follows
that 
\begin{align*}
0 &\geq \underset{ T\rightarrow \infty\,  
L\rightarrow -\infty} {\limsup}\int_M\phi g\nabla \langle\phi, \nabla g\rangle\\
&\geq\underset{ T\rightarrow \infty\,  
L\rightarrow -\infty} {\limsup}\left( \frac{n-1}{2T}\left(
\sqrt{\frac{4}{(n-1)^2}\mathcal{A}-\delta}-\varepsilon \right) \frac T2-\frac{n-1}{-2L}\left( \sqrt{\frac{4}{(n-1)^2}\mathcal{B}+\delta}+\varepsilon \right) 
\frac{-L}2\right)\\
&=\frac{n-1}4\left( \sqrt{\frac{4}{(n-1)^2}\mathcal{A}-\delta}-\sqrt{\frac{4}{(n-1)^2}\mathcal{B}+\delta}-2\varepsilon \right).
\end{align*}
Since $\varepsilon $ is arbitrary, we have \[\frac{4}{(n-1)^2}\mathcal{B}+\delta\geq \frac{4}{(n-1)^2}\mathcal{A}-\delta\] for any $\delta>0$. This proves the Proposition.

If $\mathcal{A}=\infty$, then similar argument as above will imply that $\mathcal{B}=\infty$. Then the Proposition still holds.

\begin{appendix}
\end{appendix}

%
%
%
%
%


\begin{bibdiv}
\begin{biblist}
\bib{CY75}{article}{
   author={Cheng, S. Y.},
   author={Yau, S. T.},
   title={Differential equations on Riemannian manifolds and their geometric
   applications},
   journal={Comm. Pure Appl. Math.},
   volume={28},
   date={1975},
   number={3},
   pages={333--354},
   issn={0010-3640},
   review={\MR{385749}},
   doi={10.1002/cpa.3160280303},
}
\bib{LT}{article}{
   author={Li, Peter},
   author={Tam, Luen-Fai},
   title={Harmonic functions and the structure of complete manifolds},
   journal={J. Differential Geom.},
   volume={35},
   date={1992},
   number={2},
   pages={359--383},
   issn={0022-040X},
   review={\MR{1158340}},
}

\bib{LW1}{article}{
   author={Li, Peter},
   author={Wang, Jiaping},
   title={Complete manifolds with positive spectrum},
   journal={J. Differential Geom.},
   volume={58},
   date={2001},
   number={3},
   pages={501--534},
   issn={0022-040X},
   review={\MR{1906784}},
}
\bib{LW2}{article}{
   author={Li, Peter},
   author={Wang, Jiaping},
   title={Complete manifolds with positive spectrum. II},
   journal={J. Differential Geom.},
   volume={62},
   date={2002},
   number={1},
   pages={143--162},
   issn={0022-040X},
   review={\MR{1987380}},
}
\bib{LW3}{article}{
   author={Li, Peter},
   author={Wang, Jiaping},
   title={Weighted Poincar\'{e} inequality and rigidity of complete manifolds},
   language={English, with English and French summaries},
   journal={Ann. Sci. \'{E}cole Norm. Sup. (4)},
   volume={39},
   date={2006},
   number={6},
   pages={921--982},
   issn={0012-9593},
   review={\MR{2316978}},
   doi={10.1016/j.ansens.2006.11.001},
}


\bib{V}{article}{
   author={Varopoulos, N. T.},
   title={Potential theory and diffusion on Riemannian manifolds},
   conference={
      title={Conference on harmonic analysis in honor of Antoni Zygmund,
      Vol. I, II},
      address={Chicago, Ill.},
      date={1981},
   },
   book={
      series={Wadsworth Math. Ser.},
      publisher={Wadsworth, Belmont, CA},
   },
   date={1983},
   pages={821--837},
   review={\MR{730112}},
}

\end{biblist}
\end{bibdiv}
\end{document}